
%

\magnification=1200
\hsize=11.25cm
\vsize=18cm
\parskip 0pt
\parindent=12pt
\voffset=1cm
\hoffset=1cm



\catcode'32=9

\font\tenpc=cmcsc10
\font\eightpc=cmcsc8
\font\eightrm=cmr8
\font\eighti=cmmi8
\font\eightsy=cmsy8
\font\eightbf=cmbx8
\font\eighttt=cmtt8
\font\eightit=cmti8
\font\eightsl=cmsl8
\font\sixrm=cmr6
\font\sixi=cmmi6
\font\sixsy=cmsy6
\font\sixbf=cmbx6

\skewchar\eighti='177 \skewchar\sixi='177
\skewchar\eightsy='60 \skewchar\sixsy='60

\catcode`@=11

\def\tenpoint{%
  \textfont0=\tenrm \scriptfont0=\sevenrm \scriptscriptfont0=\fiverm
  \def\rm{\fam\z@\tenrm}%
  \textfont1=\teni \scriptfont1=\seveni \scriptscriptfont1=\fivei
  \def\oldstyle{\fam\@ne\teni}%
  \textfont2=\tensy \scriptfont2=\sevensy \scriptscriptfont2=\fivesy
  \textfont\itfam=\tenit
  \def\it{\fam\itfam\tenit}%
  \textfont\slfam=\tensl
  \def\sl{\fam\slfam\tensl}%
  \textfont\bffam=\tenbf \scriptfont\bffam=\sevenbf
  \scriptscriptfont\bffam=\fivebf
  \def\bf{\fam\bffam\tenbf}%
  \textfont\ttfam=\tentt
  \def\tt{\fam\ttfam\tentt}%
  \abovedisplayskip=12pt plus 3pt minus 9pt
  \abovedisplayshortskip=0pt plus 3pt
  \belowdisplayskip=12pt plus 3pt minus 9pt
  \belowdisplayshortskip=7pt plus 3pt minus 4pt
  \smallskipamount=3pt plus 1pt minus 1pt
  \medskipamount=6pt plus 2pt minus 2pt
  \bigskipamount=12pt plus 4pt minus 4pt
  \normalbaselineskip=12pt
  \setbox\strutbox=\hbox{\vrule height8.5pt depth3.5pt width0pt}%
  \let\bigf@ntpc=\tenrm \let\smallf@ntpc=\sevenrm
  \let\petcap=\tenpc
  \normalbaselines\rm}

\def\eightpoint{%
  \textfont0=\eightrm \scriptfont0=\sixrm \scriptscriptfont0=\fiverm
  \def\rm{\fam\z@\eightrm}%
  \textfont1=\eighti \scriptfont1=\sixi \scriptscriptfont1=\fivei
  \def\oldstyle{\fam\@ne\eighti}%
  \textfont2=\eightsy \scriptfont2=\sixsy \scriptscriptfont2=\fivesy
  \textfont\itfam=\eightit
  \def\it{\fam\itfam\eightit}%
  \textfont\slfam=\eightsl
  \def\sl{\fam\slfam\eightsl}%
  \textfont\bffam=\eightbf \scriptfont\bffam=\sixbf
  \scriptscriptfont\bffam=\fivebf
  \def\bf{\fam\bffam\eightbf}%
  \textfont\ttfam=\eighttt
  \def\tt{\fam\ttfam\eighttt}%
  \abovedisplayskip=9pt plus 2pt minus 6pt
  \abovedisplayshortskip=0pt plus 2pt
  \belowdisplayskip=9pt plus 2pt minus 6pt
  \belowdisplayshortskip=5pt plus 2pt minus 3pt
  \smallskipamount=2pt plus 1pt minus 1pt
  \medskipamount=4pt plus 2pt minus 1pt
  \bigskipamount=9pt plus 3pt minus 3pt
  \normalbaselineskip=9pt
  \setbox\strutbox=\hbox{\vrule height7pt depth2pt width0pt}%
  \let\bigf@ntpc=\eightrm \let\smallf@ntpc=\sixrm
  \let\petcap=\eightpc
  \normalbaselines\rm}
\catcode`@=12

\tenpoint



\catcode`\@=11
\def\pc#1#2|{{\bigf@ntpc #1\penalty \@MM\hskip\z@skip\smallf@ntpc%
	\uppercase{#2}}}
\catcode`\@=12

\def\pointir{\discretionary{.}{}{.\kern.35em---\kern.7em}\nobreak
   \hskip 0em plus .3em minus .4em }

\def\qed{\quad\raise -2pt\hbox{\vrule\vbox to 10pt{\hrule width 4pt
   \vfill\hrule}\vrule}}

\def\rem#1|{\par\medskip{{\it #1}\pointir}}

\def\vspace[#1]{\noalign{\vskip#1}}

\def\abstract#1{\vbox{\eightpoint\narrower\narrower 
\pc ABSTRACT|\pointir #1}}


\def\setN{{\enslettre N}}

\long\def\maskbegin#1\maskend{}


\def\section#1{\goodbreak\par\vskip .3cm\centerline{\bf #1}
   \par\nobreak\vskip 3pt }

\long\def\th#1|#2\endth{\par\medbreak
   {\petcap #1\pointir}{\it #2}\par\medbreak}

\def\article#1|#2|#3|#4|#5|#6|#7|
    {{\leftskip=7mm\noindent
     \hangindent=2mm\hangafter=1
     \llap{{\tt [#1]}\hskip.35em}{\petcap#2}\pointir
     #3, {\sl #4}, {\bf #5} ({\oldstyle #6}),
     pp.\nobreak\ #7.\par}}
\def\livre#1|#2|#3|#4|
    {{\leftskip=7mm\noindent
    \hangindent=2mm\hangafter=1
    \llap{{\tt [#1]}\hskip.35em}{\petcap#2}\pointir
    {\sl #3}, #4.\par}}
\def\divers#1|#2|#3|
    {{\leftskip=7mm\noindent
    \hangindent=2mm\hangafter=1
     \llap{{\tt [#1]}\hskip.35em}{\petcap#2}\pointir
     #3.\par}}


\def\frac#1#2{{#1\over #2}}

\def\JFracSym{\bf J}
\def\JFrac#1#2{{\JFracSym}\Bigl[ \matrix{ #1\cr #2\cr } \Bigr] }
\def\bfu{{\bf u}}
\def\bfv{{\bf v}}
\def\bfa{{\bf a}}
\def\bfJ{{\bf J}}

\def\setN{{\bf N}}

\font\KFracFont=cmr12 at 20pt
\def\KK{\mathop{\lower 4pt\hbox{\KFracFont K}}\limits}

\long\def\maskbegin#1\maskend{}

%
\def\ThmJFGene{3.3}

\def\ThmJFGeneLevel{3.4}
\def\ThmJFGeneLevelK{3.5}
\def\EqJFQuadra{3.1}
\def\EqJFChangeVar{3.2}
\def\EqJFChangeVarK{3.3}

\rightline{May 16, 2014}
\bigskip


\centerline{\bf Hankel Determinant Calculus }
\centerline{\bf for the Thue-Morse and related sequences}

\bigskip
\centerline{\sl Guo-Niu Han}
\footnote{}{\eightpoint
}
\bigskip
\bigskip
{\narrower\narrower
\eightpoint
\noindent
{\bf Abstract}.\quad
The Hankel determinants of certain automatic sequences $f$ are evaluated, 
based on a calculation modulo a prime number.
In most cases, the Hankel determinants of automatic sequences
do not have any closed-form expressions; the traditional methods,
such as $LU$-decompo\-si\-tion and Jacobi continued fraction, cannot be applied directly.
Our method is based on a simple idea:
the Hankel determinants of 
each sequence $g$ equal to $f$ modulo $p$  are equal to the Hankel
determinants of $f$ modulo $p$. 
The clue then consists of finding a nice 
sequence $g$, whose Hankel determinants have closed-form expressions.

Several examples are presented, including a result  saying that the Hankel determinants
of the Thue-Morse sequence are nonzero, 
first proved by
Allouche, Peyri\`ere, Wen and Wen using determinant manipulation.
The present approach shortens 
the proof of the latter result significantly.
We also prove that the corresponding 
Hankel determinants do not vanish
when the 
powers $2^n$ in the infinite product defining the $\pm 1$ Thue--Morse 
sequence are replaced by $3^n$.

}

\section{1. Introduction} 

Let $x$ be a parameter. We identify a sequence ${\bf a}=(a_0, a_1, a_2, \ldots)$ and its generating function 
$f=f(x)=a_0+a_1x+a_2x^2+\cdots$. Usually, $a_0=1$.
For each $n\geq 1$ and $k\geq 0$
the Hankel determinant of the series $f$ (or of the sequence $\bf a$)  is defined by
$$
H_n^{(k)} (f) := \left|
\matrix{ a_k & a_{k+1} & \ldots & a_{k+n-1} \cr
a_{k+1} & a_{k+2} & \ldots & a_{k+n} \cr
\ \vdots \hfill & \ \vdots \hfill & \ddots &
\ \vdots \hfill \cr
a_{k+n-1} & a_{k+n} & \ldots & a_{k+2n-2} \cr} \right|.
\leqno{(1.1)}
$$
Let $H_n(f):=H_n^{(0)}(f)$, for short;
the {\it sequence of the Hankel determinants} of $f$ is defined to be: 
$$H(f):=(H_0(f)=1, H_1(f), H_2(f), H_3(f), \ldots).$$ 
In some cases
Hankel determinants can be evaluated by using 
basic determinant manipulation, $LU$-decomposition, or  Jacobi 
continued fraction (see, e.g., [Kr98, Kr05, Fl80, Wa48, Mu23]).
However, the Hankel determinants of several power series $f$ 
related to automatic sequences
do not seem to have closed-form expressions, 
as will be seen in this paper. 
The following result by
Allouche, Peyri\`ere, Wen and Wen [APWW] in 1998,
has
strongly motivated the present paper.

\proclaim Theorem 1.1 [APWW].
Let $P_2=P_2(x)=\prod_{k=0}^\infty(1-x^{2^k})$ be the $\pm 1$ Thue-Morse sequence. Then
$H_n(P_2)\not=0$ for every positive integer $n$.

The first values of the coefficients and Hankel determinants of $P_2(x)$
are: 
$$\leqalignno{
	P_2&=(1, -1, -1, 1, -1, 1, 1, -1, -1, 1, 1, -1, 1, -1, -1, 1, -1, \ldots)\cr
	H(P_2)&= (1, 1, -2, 4, 8, -16, -32, -64, 128, -256, -1536, -3072,\ldots)\cr
}$$

A combinatorial proof of Theorem 1.1 was recently derived by Bugeaud and the author [BH13]. 

\smallskip

Let ${\bf u}=(u_1, u_2, \ldots)$ and ${\bf v}=(v_0, v_1, v_2, \ldots)$
be two sequences.
Recall that the {\it Jacobi continued fraction} attached to $(\bfu, \bfv)$, 
or {\it $J$-fraction}, for short, is a continued
fraction of the form
$$
{v_0 \over 1 + u_1 x - 
	\displaystyle{v_1 x^2 \over 1+u_2x - 
		\displaystyle{v_2 x^2 \over {1 + u_3x - 
\displaystyle{v_3 x^2\over \ddots}}} }}, 
\leqno{(1.2)}
$$
also denoted by
$$
\JFrac{\bf u}{\bf v}=
{\bf J}[{{\bf u} / {\bf v}}]=
\JFrac{u_1, u_2, \cdots}{v_0, v_1, v_2, \cdots}.
$$
The basic properties on $J$-fractions, we now recall, can be found in
[Fl80, Wa48, Vi83].
The $J$-fraction of a given power series  $f$ exists (i.e., 
$f= \bfJ[\bfu / \bfv]$) if and only if
all the Hankel determinants $H_n(f)$ of $f$ are nonzero.
The first values of the coefficients $u_n$ and $v_n$ in the 
$J$-fraction expansion can be calculated by
the {\it Stieltjes Algorithm}.
Also, Hankel determinants can be calculated
from the $J$-fraction by means of the following  fundamental relation: 
$$
H_n\Bigl(\JFrac{u_1, u_2, \cdots }{ v_0, v_1, v_2, \cdots}\Bigr)
= v_0^n v_1^{n-1} v_2^{n-2} \cdots v_{n-2}^2 v_{n-1}. \leqno{(1.3)}
$$
Conversely, the coefficients $u_n$ and $v_n$ in the $J$-fraction can be 
calculated using the Hankel determinants by means of the  following relations, 
when all denominators are nonzero. 
$$
\leqalignno{
	u_{n} &= -{1\over H^{(1)}_{n-1}} \Bigl(
	{ H_{n-1}  H^{(1)}_{n}\over H_{n}}
	+ { H_{n}  H^{(1)}_{n-2}\over H_{n-1}}
	\Bigr), \qquad (n\geq 2) &{(1.4)}\cr
	\noalign{\smallskip}
	v_{n} &= 
	{ H_{n}  H_{n-2}\over (H_{n-1})^2}. \qquad (n\geq 2) &(1.5)\cr
}
$$
Relation (1.3) is an efficient method for evaluating Hankel determinants. 
\goodbreak

Let us try to evaluate 
the Hankel determinants for the Thue-Morse sequence
by using the $J$-fraction. By the Stieltjes algorithm, we get  
$$
P_2(x)=
\JFrac{\bf u}{\bf v}=
\JFrac
{ 1, -1, 1, -1, 1, -1, 1, -1, 1, -1, 1, -1 \cdots}
{1, -2, 1, -1, -1, -1, 1, -1, 1, -3, {1\over 3}, -{1\over 3}, -3 \cdots}.
$$
The top coefficients $u_n$ seem to be  very simple.
However, we are not able to guess any closed-form expression for the bottom 
coefficients $v_n$, which are even rational numbers. Therefore, we cannot prove anything 
about the Hankel determinants.
\medskip

Coons [Co13], using the method described in [APWW],
proved the following theorem.

\proclaim Theorem 1.2 [Coons].
Let
$$
S_2=S_2(x)= {1\over x} \sum_{n=0}^\infty {x^{2^n} \over 1+x^{2^n}}.
$$
Then $H_n(S_2)\equiv 1\pmod 2$.

Again, we are not able to guess any closed-form expression for the Hankel
determinants of $S_2$, as the 
first values of the coefficients of the series,  the Hankel determinants and the
$J$-fraction of $S_2$ read: 
$$
\leqalignno{
	S_2&=(1, 2, 1, 3, 1, 2, 1, 4, 1, 2, 1, 3, 1, 2, 1, 5, 1, 2, 1, 3, 1, 2, 1, 4, \ldots) \cr
	H(S_2)&=(1, 1, -3, -1, 21, 1, -3, -9, 945, 9, -3, -1, 21, 9, -243,  \ldots)\cr
	S_2&=
\JFrac
{ 
-2, {7\over 3}, {23\over 3}, -{167\over 21}, -{169\over 21}, 7, 7, 
-{629\over105}, -{631\over 105}, 7, 7, -{57\over7}, -{55\over 7}, 
\cdots}
{
1,	-3, -{1\over 9}, -63, -{1\over 441}, -63, -1, -35, -{1\over 11025}, 
	-35, -1, -63, \cdots}
\cr
}
$$
The main idea to solve the problem is to proceed  as
follows:
\medskip

let $p$ be a prime number and $f$ a sequence. We want to prove that
$H_n(f)\not= 0 \pmod p$; 
if, apparently, there is no closed-form for the coefficients in the $J$-fraction of $f$, 
we try to find a sequence $g\equiv f \pmod p$, 
such that the Hankel determinants of $g$ have a closed form.
As it is easy to prove that
$H_n(f)\equiv H_n(g)\pmod p$,
it is very likely that some properties on the Hankel determinants of $f$ 
can be established.

\medskip
\proclaim Question. 
{\it How to find a {\it nice} sequence $g$ such that $g\equiv f$ 
for which each coefficient in the $J$-fraction of $g$ has a closed-form expression?} 

\medskip
By observing the occurrences of the factor $2$ in the coefficients of the $J$-fraction of $S_2$ given in Theorem 1.2,
we guess the following 
``nice" sequence 
$$
g=\JFrac{0,1,1,1,1,1,\ldots}{1,1,1,1,1,1,1\ldots},
$$
whose Hankel determinant is $H_n(g)=1\not=0$. For proving Theorem~1.2, it remains to prove that
$S_2 \equiv g \pmod 2$. 
For Theorem~1.1, it is more complicated; 
we need the so-called {\it grafting technique}. 
The proofs of Theorems 1.1-2 are given in Section 2
with further examples.
In Section 3 we derive two $J$-fractions by using the chopping method
(Proposition 3.2$''$ and Theorem 3.3) and prove that 
the Hankel determinant sequences of several power series are periodic
(Propositions 3.6-8). 

On the one hand, we provide short proofs of results established 
in the papers [APWW, Co13],  on the other hand, we obtain several new results.
In particular, we should like to single out the following theorem.

\proclaim Theorem 1.3. 
Let $P_3=P_3(x)=\prod_{k\geq 0} (1-x^{3^k})$.
Then $H_{n}(P_3) \equiv (-1)^{n} \pmod 3$ for every positive integer $n$.

Notice that the sequence $P_3$ is obtained from 
the Thue-Morse sequence $P_2$ by modifying the 
exponent of $x$ from 2 to 3. 
It is worth mentioning that,
when $m\geq 4$,
the Hankel determinants for the sequence $\prod_{k\geq 0} (1-x^{m^k})$ 
are not all nonzero. 
A self-contained and short proof of Theorem~1.3 is found in Section~4. 
The following result, that could be called ``{\it one sequence, two modulos}",
is also proved in Section~4. 
\proclaim Theorem 1.4.
We have
$$
\leqalignno{
	\sqrt{1\over (1-x)(1+3x)} &\equiv \prod_{k=0}^\infty (1-x^{2^k})  \pmod 4,\cr
	\sqrt{1\over (1-x)(1+3x)} &\equiv \prod_{k=0}^\infty (1-x^{3^k})  \pmod 3.\cr
}
$$


\section{2. Hankel determinants modulo $p$ and the grafting technique} 

Let $p$ be a prime number.
For a given power series $f$ we present some methods for guessing {\it and} 
calculating the $J$-fraction 
of $f$,
and also proving properties mod $p$ for its Hankel determinants. 
An ultimately periodic sequence is written in contracted form by using the star sign. For instance, the sequence $\bfa=(1, (3,0)^*)$ 
represents $(1,3,0,3,0,3,0,\ldots)$, that is, $a_0=1$ and $a_{2k+1}=3,
a_{2k+2}=0$ for each positive integer $k$. Two sequences $\bf a$ and $\bf b$ are said to be 
{\it congruent modulo $p$} if $a_k\equiv b_k \pmod p$ for all $k$. 
For each integer $z$ we have $(x+z)^p\equiv x^p+z^p\pmod p$ and derive
the following lemma.

\proclaim Lemma 2.1. Let $f(x)$ be a power series with integral coefficients. Then
$$f(x)^p\equiv f(x^p)\pmod p.\leqno{(2.1)}$$

Let $a_1,b_1, a_2, b_2$ be four integers such that 
$(p,b_1)=1$ and  $(p,b_2)=1$. 
The two fractions ${a_1/b_1}$ and ${a_2/b_2}$
are said to be {\it congruent modulo} $p$ if
$a_1b_2\equiv a_2b_1\pmod p$. We write
${a_1/b_1}\equiv {a_2/b_2} \pmod p$.
This fractional congruence is closed under addition and multiplication.
Let $a_1/b_1\equiv c_1\pmod p$ and $a_2/b_2\equiv c_2\pmod p$, then
$a_1/b_1+a_2/b_2 \equiv c_1+c_2 \pmod p$
and
$a_1/b_1\times a_2/b_2 \equiv c_1c_2 \pmod p$.
The fractional congruence for power series is also closed under addition and multiplication. Also, the ring of formal power series with rational coefficients
modulo $p$ 
is an integral domain. 

\proclaim Lemma 2.2. Let $f$ and $\hat f$ be two power series with rational coefficients and
$\bfJ[\bfu, \bfv]=f,\, \bfJ[\hat\bfu, \hat\bfv]=\hat f$ be their $J$-fraction expansions. 
Then
\smallskip
(1) If $f\equiv \hat f\pmod p$,  then $H(f)\equiv H(\hat f) \pmod p$.
\smallskip
(2) If $\bfu\equiv \hat \bfu\pmod p$ and $\bfv\equiv \hat \bfv\pmod p$,  
then $f\equiv \hat f \pmod p$.
\smallskip
(3) If $\bfv\equiv \hat\bfv\pmod p$,  then $H(f)\equiv H(\hat f) \pmod p$.

{\it Proof}. (1) The Hankel determinants are expressed in terms of the coefficients of the power series by using only addition and multiplication.
(2)~The coefficients of the power series are expressed in terms of the coefficients in the $J$-fraction by using only addition and multiplication.
(3)~By the fundamental relation (1.3).
\qed
\smallskip

{\it Remark}. The converse of (1) is not true.  A counter-example is the following pair with $p=2$:
$$
f={1-\sqrt{1-{4x^2\over 1-x}}\over 2x^2}
\hbox{\quad and\quad }
\hat f={1-\sqrt{1-{4x\over 1-x}}\over 2x}.
$$
\medskip

Let $f$ be a power series and $g$ be a $J$-fraction.
If the two sequences $\bfu$ and $\bfv$ in 
$g$ are ultimately periodic with the same period,
we can check that $f$ and $g$ are equal or not.
For example,
we claim that the $J$-fraction of the power series
$$
f={(1-x)(1+2x)-\sqrt{(1-x)(1-2x)(1+3x)(1+2x-4x^2)} \over 4x^2(1-x)}
$$
is equal to
$$
g=\JFrac
{ \noalign{\smallskip} (-{1\over 2}, -{1\over 2}, 2)^*}
{1, ({1\over 4}, 2,2)^*}
.
$$
To see this,
we check that  $f$ verifies the following quadratic functional equation
$$
f={1\over 1- {1\over 2}x-
	\displaystyle{{1\over 4}x^2\over 1- {1\over 2}x- 
\displaystyle{2x^2\over 1+ 2x-2x^2f}}}.
$$
Moreover, the first values of $f$ and $g$ are the same,  namely,
$(1, {1\over 2}, {1\over 2}, \ldots)$.
Hence,
the two power series $f$ and $g$ are equal.
Later in the paper 
this kind of proof will not be reproduced,
as it can be
done automatically: the sentence {\it ``we can prove"} replaces the full proof.

\medskip

{\it Proof of Theorem 1.2 [Coons]}. We have
$$(xS_2(x))^2\equiv x^2S_2(x^2) = 
\sum_{n=1}^\infty {x^{2^n} \over 1+x^{2^n}}=xS_2(x)-{x\over 1+x}\pmod 2$$
so that
$$xS_2(x)^2\equiv S_2(x)-{1\over 1+x}\pmod 2$$
and
$$ \Bigl(S_2(x)-{1+\sqrt{1-3x\over 1+x}\over 2x}\Bigr) 
 \Bigl(S_2(x)-{1-\sqrt{1-3x\over 1+x}\over 2x}\Bigr)\equiv 0 \pmod 2.$$
We get 
$$ S_2(x)\equiv{1-\sqrt{1-3x\over 1+x}\over 2x}\pmod 2. $$
Let $g$ be the right-hand side of the above equation. We can prove 
$$
g=
\JFrac
{ 0, (-1)^* }
{ (1)^*}
.
$$
Hence, $H_n(g)\equiv 1\pmod 2$, so does $H_n(S_2)$ by Lemma 2.2(1). \qed
\medskip

{\it Proof of Theorem 1.3}.  We successively have 
$$\displaylines{
	P_3(x)=(1-x)P_3(x^3)\equiv (1-x)P_3(x)^3 \pmod 3,\cr
P_3(x) ( 1- (1-x)P_3(x)^2) \equiv 0\pmod 3,\cr
 1- (1-x)P_3(x)^2 \equiv 0\pmod 3,\cr
P_3(x)^2 \equiv {1\over 1-x }\pmod 3,\cr
\Bigl(P_3(x) - \sqrt{1\over 1-x }\Bigr) \Bigl(P_3(x) + \sqrt{1\over 1-x }\Bigr) \equiv 0\pmod 3,\cr 
P_3(x) \equiv \sqrt{1\over 1-x} \pmod 3.\cr 
}
$$
Notice that $P_3(x)$ has integral coefficients, but 
$\sqrt{{1\over 1-x}}$ has rational coefficients.
We can prove that
$$
\sqrt{1\over 1-x}=\JFrac
{ (-1/2)^*}
{1, 1/8, (1/16)^*}
.
$$
The above $J$-fraction itself
is congruent to
$$
g=
\JFrac
{ (1)^*} 
{1, -1, (1)^*}
$$
modulo 3, by Lemma 2.2(2), knowing that $1/2 \equiv -1 \pmod 3$. We have 
$H_n(g)=(1, -1)^*$. Hence, $H_n(P_3)\equiv H_n(g)\equiv (1,-1)^* \pmod 3$. \qed
\smallskip
There is also a proof without using fractional congruence. See Section~4.
Notice that $H(P_3(x))=H(P_3(-x))$ by (1.3). 
That means $H(g)\not=0$ for 
$g=\prod_{k\geq 0} (1+x^{3^k})$. 

\medskip

For proving Theorem 1.1, we need a technique, called ``{\it grafting}". 
Let $F(x)$ and $G(x)$ be two $J$-fractions
$$
F(x)=\JFrac{u_1, u_2, u_3, \cdots}{v_0, v_1, v_2, v_3, \cdots}
\hbox{\quad \rm and\quad}
G(x)=\JFrac{a_1, a_2, a_3, \cdots}{b_0, b_1, b_2, b_3, \cdots}
$$
such that $b_0=1$. For each $k\in\setN$ the {\it grafting} of $G(x)$ into $F(x)$
of order~$k$, denoted by $F(x) |^r G(x)$, is defined to be the following 
$J$-fraction
$$
F(x) |^k G(x)=
\JFrac{u_1, u_2,  \cdots, u_{k}, a_1, a_2, a_3, \cdots}
{v_0, v_1, v_2, \cdots, v_{k}, b_1, b_2, b_3, \cdots}.
$$
Let $F|G:=F|^1G$ and $F||G:=F|^2G$, for short.

If $u_i, v_i \pmod p$ exists and $v_i\not\equiv 0\pmod p$ for all $i\geq k+1$, we define
$$\bar G:= \JFrac{u_{k+1} \pmod p, u_{k+2}\pmod p, u_{k+3}\pmod p, \cdots}
{1, v_{k+1}\pmod p, v_{k+2}\pmod p, v_{k+3}\pmod p, \cdots}
$$
and $\bar F=F|^r \bar G$. Then the Hankel determinants of $F$ and $\bar F$
have the following relation 
$${H_n(F)\over H_n(\bar F)} \equiv 1 \pmod p\leqno{(2.2)}$$
in view of the fundamental relation (1.3).

For instance, 
the first values of the $J$-fraction of the Thue-Morse sequence $P_2$ are
$$
P_2=\JFrac
{ (1, -1)^* }
{1, -2, 1, -1, -1, -1, 1, -1, 1, -3, {1\over 3}, -{1\over 3}, -3, 
	1, -1, 1, 1, -3,
\cdots}
.
$$
We see that the previous sequences $\bfu, \bfv$ contain only one even number,
$-2$, and it occurs 
at position $v_1$.
Delete $(v_1, u_1)$, which means that 
we define
the following $J$-fraction $g$
$$
g=\JFrac
{ (-1, 1)^* }
{1, 1, -1, -1, -1, 1, -1, 1, -3, {1\over 3}, -{1\over 3}, -3,
 1, -1, 1, 1, -3,
\cdots}
, 
$$
so that all the Hankel determinants of $g$ are odd fractional numbers by~(1.3).
\medskip

{\it Proof of Theorem 1.1}.
Define the sequence $g$ by
$$
P_2={1\over 1+x+2x^2g},
$$
or
$$
g={1\over 2x^2}({1\over P_2} -1-x).
$$
By Theorem 1.4 the following identities hold: 
$$
\displaylines{
1/P_2 \equiv \sqrt{(1-x)(1+3x)} \pmod 4,\cr
g\equiv {1\over 2x^2}(1+x-\sqrt{(1-x)(1+3x)})\pmod 2.\cr
}
$$
We can prove that the right-hand side $\bar g$ of the above equation has a simple $J$-fraction
$$
\bar g=\JFrac
{ (1)^* }
{ (1)^* }.
$$
Let $\bar P_2$ be the grafting of $\bar g$ into $P_2$
$$\bar P_2=P_2|\bar g=     
{1\over 1+x+2x^2\bar g},
$$
so that
$H_n(\bar P_2)=(-2)^{n-1}$ from (1.3).
Hence, $H_n(P_2)/2^{n-1}\equiv 1 \pmod 2$ by (2.2).\qed
\medskip
Let $P_2=\sum_{n=0}^\infty \eta_n x^n$ be the Thue-Morse sequence. We now evaluate
the Hankel determinants of the following two sequences
$$
\leqalignno{
	\delta_n & =(\eta_{n} - \eta_{n+1})/2, &(2.3) \cr
	\gamma_n & =(\eta_{n} - \eta_{n+2})/2. &(2.4) \cr
}
$$

\smallskip
The following result was proved  in [APWW, Proposition 2.2(2)].

\proclaim Proposition 2.3.
The Hankel determinants of the sequence $(\delta_n)_{n=0,1,2,\ldots}$ 
are odd integral numbers.

{\it Proof.} 
The generating function for the sequence $(\gamma_n)$ is equal to
$$
f={1-(1-x)P_2\over 2x},
$$
which is congruent to
$$
g:={1-(1-x){\sqrt{1\over (1-x)(1+3x)}}\over 2x} \pmod 2.
$$
by Theorem 1.4.
We can prove that $g$ has the following $J$-fraction expansion
$$
g={{1-\sqrt{1-x\over 1+3x}}\over 2x} 
=\JFrac
{ 2, (1)^* }
{ (1)^* }
.
$$
Hence, $H_n(g)=1$ and $H_n(f)\equiv 1 \pmod 2$. \qed

\proclaim Proposition 2.4.
The Hankel determinants of the sequence $(\gamma_n)_{n=0,1,2,\ldots}$
are odd integral numbers.

{\it Proof.}
The generating function for the sequence $(\gamma_n)$ is equal to
$$
f={1-x-(1-x^2)P_2\over 2x^2},
$$
which is congruent to
$$
g:=-{1-x-(1-x^2){\sqrt{1\over (1-x)(1+3x)}}\over 2x^2} \pmod 2
$$
by Theorem 1.4.
We can prove that $g$ has the following $J$-fraction expansion
$$
g=-{{1-x-(1+x)\sqrt{1-x\over 1+3x}}\over 2x^2} 
=\JFrac
{ (3,-1)^* }
{1, (-1)^* }
.
$$
Hence, $H_n(g)\equiv 1\pmod 2$ and $H_n(f)\equiv 1 \pmod 2$. \qed

\proclaim Proposition 2.5.
Let 
$$f=3\prod_{n=1}^\infty(1-x^{3^n})-{2\over 1-x}.$$
Then, $H_k(f)\not=0$ for all $k$.

{\it Remark}. When replacing the factor $1-x^{3^n}$ by $1+x^{3^n}$ in the above formula, experimental calculation of the first values suggests that
all the Hankel determinants are still nonzero.
However, we are not able to prove 
that the latter Hankel determinants do not vanish.
\medskip  

{\it Proof}.
We have
$$
f=\JFrac
{ 2, -7/2, 7/10, 32/65, -187/26, 259/34, -49/272, 241/16,\ldots }
{1, -6, -5/4, -26/25, 10/169, -221/4, 64/289, -17/256, \ldots}
.
$$
The factor $3$ occurs only once, at position $v_1$.
We use the grafting technique.
Define
$$
f={1\over 1+2x+6x^2g},\leqno{(2.5)}
$$
or
$$
(1+2x+6x^2g)\bigl(3\prod_{n=0}^\infty(1-x^{3^n})-{2}\bigr)= 1-x.\leqno{(2.6)}
$$
By (1.3) we have
$$H_n(f)=(-6)^{n-1}H_{n-1}(g).\leqno{(2.7)}$$
for all $n$.
By Theorem 1.4 identity  (2.6) becomes
$$
(1+2x+6x^2g)\Bigl(3 \sqrt{1\over (1-x)(1+3x)} -{2}\Bigr)\equiv 1-x,\pmod 9
$$
or
$$
g\equiv 
{(1+2x)\sqrt{1\over (1-x)(1+3x)} -1-x \over x^2}
\pmod 3.
$$
Let $h$ be the right-hand side of the above equation. 
Then,
$$
h=\JFrac
{ (-1)^* }
{ 1, (4, -1/2, -1/2)^*}
,
$$
so that $H_n(h)\equiv 1\pmod 3$.
Hence, $H_n(g)\equiv 1\pmod 3$. By (2.7) we have $H_n(f)\not=0$. \qed

\proclaim Proposition 2.6.
Let $f$ be the sequence obtained from $P_3$ by deleting the first term, i.e.,
$f=(1-P_3)/x$.
Then, $H(f)\equiv (1)^* \pmod 3$.

{\it Proof}.
By Theorem 1.4,  
$$
f\equiv {1-\sqrt{1\over (1-x)(1+3x)}\over x} \pmod 3.
$$
Let $g$ be right-hand side of the above equation. 
Then,
$g$ has the $J$-fraction 
$$
g=\JFrac
{ 3, (-2, 5/2, 5/2)^* }
{1, (-2,-2, 1/4)^*}
	\equiv
\JFrac
{ 0, (-1, 1, 1)^* }
{(1)^*}
	\pmod3.
$$
Hence, $H(g)\equiv (1)^* \pmod 3$, so does $H(f)$.\qed

\proclaim Theorem 2.7.
Let
$$
f=f(x)=\prod_{k\geq 0}(1-x^{3^k} - x^{2\cdot3^k}).\leqno{(2.8)}
$$
Then
$H_n(f)\not=0$. 

{\it Proof}. We successively have 
$$
\displaylines{
f(x^3)=\prod_{k\geq 1}(1-x^{3^{k}} - x^{2\cdot3^{k}});\cr
f=f(x)=\prod_{k\geq 0}(1-x^{3^k} - x^{2\cdot3^k})
=(1-x-x^2) f(x^3);\cr
f\equiv (1-x-x^2) f^3 \pmod 3; \cr
1\equiv (1-x-x^2) f^2 \pmod 3; \cr
f\equiv \sqrt{1/(1-x-x^2)}  \pmod 3. \cr
}
$$
The right-hand side of the above equality has the following $J$-fraction
expansion:
$$ \sqrt{1\over 1-x-x^2}  
=
\JFrac
{(-1/2)^*}
{1,5/8, (5/16)^*}
\equiv 
\JFrac
{(1)^*} 
{1, 1, (-1)^*}
\pmod 3,
$$
so that
$$
H(f) \equiv (1,(1,1,2,2)^*) \pmod 3.\qed\leqno{(2.9)}
$$

{\it Remark}. The sequence $f$ defined in (2.8) is a $\{1, -1\}$-sequence.


\section{3. Continued fraction and the chopping method} 
\medskip
When the two coefficients in the $J$-fraction are
ultimately periodic with the same period, the corresponding power series
is easy to obtained. However this is not always the case, as shown in the following Proposition. 

\proclaim Proposition 3.1.
Let
$$
f(x)={1-\sqrt{1-{4x^4\over 1-x^2}}\over 2x^4}.
$$
Then
$$
f=\JFrac
{ (0)^* }
{1, 1,1,-1,-1, 2, 1/2, -1/2,-2, 3, 1/3, -1/3, -3,\ldots}.
$$
In other words, if $f=\bfJ[\bf u/\bf v]$, then
$u_k=0$ and $v_{4k+1}=k, v_{4k+2}=1/k, v_{4k+3}=-1/k, v_{4k+4}=-k$ for every positif integer $k$.

The proof of Proposition 3.1 is based on the following 
generalization with one more parameter $z$. 
Proposition 3.1$'$ becomes Proposition 3.1 when $z=1$. 
\proclaim Proposition 3.1$'$.
Let
$$
f=f(x;z)={1-(2z-1)x^2-\sqrt{(1-x^2)(1-x^2-4x^4)}\over 
2x^2((1-z)+(1-z+z^2)x^2-x^4)}.
$$
Then
$$
f=\JFrac
{ (0)^* }
{1, z,1/z,-1/z,-z, z+1, 1/(z+1), -1/(z+1),-(z+1), \ldots}
.
$$
In other words, if $f=\bfJ[\bf u/\bf v]$, then
$u_k=0$ and $v_{4k+1}=z+k, v_{4k+2}=1/(z+k), v_{4k+3}=-1/(z+k), v_{4k+4}=-(z+k)$ for every positif integer~$k$.

{\it Proof}. We need to check that $f(x;z)$ verifies the following
functional equation:
$$
f(x;z)={1\over
	1-\displaystyle{zx^2\over
		1-\displaystyle{{1\over z}x^2\over
			1+\displaystyle{{1\over z}x^2\over
				1+{zx^2 f(x; z+1)
}}}}
}.\qed
$$

Let us explain how to get Proposition 3.1$'$ from Proposition 3.1.
Let
$$
f_1={1-\sqrt{1-{4x^4\over 1-x^2}}\over 2x^4}
$$
and
$$
f_1=\JFrac
{ (0)^* }
{1, 1,1,-1,-1, 2, 1/2, -1/2,-2, 3, 1/3, -1/3, -3,\ldots}
.
$$
Define $f_2$ by deleting the first four pairs $u_i, v_i$ ($i=1,2,3,4$)
from the $J$-fraction of $f_1$. In other words,
$$
f_2=\JFrac
{ (0)^* }
{1, 2, 1/2, -1/2,-2, 3, 1/3, -1/3, -3,\ldots}
.
$$
By the very definition of the continued fraction we get the first values of~$f_2$
$$
f_2=(1, 0, 2, 0, 5, 0, 12, 0, 30, 0, 75, 0, 190, 0, 483, 0, 1235, 0, 3167,\ldots).
$$
With the help of a computer algebra system (see [Ru06] for example), we observe that
$f_2$ satisfies the equation
$$
(x^6-3x^4+x^2)f_2^2+(-3x^2+1)f_2-1=0.
$$
Define $f_3$ by deleting the first four pairs $u_i, v_i$ ($i=1,2,3,4$)
from the $J$-fraction of $f_2$ and repeat these steps,
we  sucessively get 
$$
\displaylines{
	(x^6-7x^4+2x^2)f_3^2+(-5x^2+1)f_3-1=0,\cr
	(x^6-13x^4+3x^2)f_4^2+(-7x^2+1)f_4-1=0,\cr
	\cdots\cr
}
$$
and guess the general equation valid for every $z$
$$
(x^6-(z^2-z+1)x^4+(z-1)x^2)f_z^2+(-(2z-1)x^2+1)f_z-1=0.
$$
Solving the above equation yields the series $f(x;z)$, defined in Proposition~3.1$'$.
The above procedure of finding generalization of $J$-fraction will be
called the {\it chopping method}.

\proclaim Proposition 3.2. Let
$$
f(x)={1-\sqrt{1-{4x^4\over 1+x}}\over 2x^4}.
$$
Then 
$$
H(f(x))=(1,1,0,0,-1,-1,0,0)^*
$$
\goodbreak

As $H_3(f(x))=0$, the traditional method fails. We then have to 
find a polarization, as stated in the following example, which becomes Proposition~3.2 when $y=-1$ and $z=0$. Notice that Proposition 3.1 is also a special case of Proposition 3.2$'$ by taking $y=0$ and $z=1$.

\proclaim Proposition 3.2$'$. Let
$$
f(x; y,z)={1-\sqrt{1-{4x^4\over 1-yx-zx^2}}\over 2x^4}.
$$
Then
$$
f(x;y,z)=\JFrac
{ (-y,0,0,0)^* }
{1, z,1/z,-1/z,-z, 2z, 1/(2z), -1/(2z),-(2z), \ldots}
.
$$

By Proposition 3.2$'$ and the fundamental relation (1.3), the Hankel determinants of $f(x;y,z)$ are 
$$
H(f(x;y,z))=(1,1,z,z,-1,-1,-2z,-2z,1,1,3z,3z, \ldots).
$$
When $z=0$ and $y=-1$ we get 
$$
H(f(x;-1,0))=(1,1,0,0,-1-1,0,0)^*.
$$
Proposition 3.2 is proved.
\medskip
However we are not able to prove Proposition 3.2$'$ directly.
By using the chopping method we find and prove the following generalization
of Proposition~3.2$'$.
Letting $t=0$ in Proposition 3.2$''$ we get Proposition~3.2$'$.

\proclaim Proposition 3.2$''$. Let
$$
f(x)=
-\,{\frac {2\,zt{x}^{2}+z{x}^{2}+yx-1+\sqrt {
\left( 4\,{x}^{4}+yx-1+z{x}^{2} \right)  \left( yx-1+z{x}^{2} \right)
}}{2{x}^{2} \left( -z{x}^{4}-{x}^{3}y+{x}^{2}+{x}^{2}{z}^{2}t+{x}^
{2}{z}^{2}{t}^{2}+yztx-zt \right) }}
$$
Then
$$
f=\JFrac
{ (-y,0,0,0)^* }
{1, (t+1)z,{1\over (t+1)z},-{1\over (t+1)z},-(t+1)z, (t+2)z, {1\over (t+2)z},  \ldots}
.
$$

\medskip

By using the  chopping methods, we derive the following
continued fraction.

\proclaim Theorem 3.3.
Let
{
$$
g=
{\frac {-2z{x}^{2}-(sx-{x}^{2}y-1)-\sqrt { 
(sx -x^2y-1)^2-(2x^2)^2  
}}{2{x}^{2}
( {x}^{2} +{x}^{2}{z}^{2} +z( sx -{x}^{2}y-1 )) }}.
$$
}Then
$$
g=\JFrac
{ (-s, 0)^* }
{ v_0, v_1, v_2, \ldots}
=\displaystyle{1\over 1- sx -
     \displaystyle{ 
			 {\alpha_1\over \alpha_0} x^2\over 
	 1-  
	\displaystyle{
	{\alpha_0\over \alpha_1}x^2
		\over 1- sx -
     \displaystyle{ 
			 {\alpha_2\over \alpha_1} x^2\over 
			 1-\displaystyle{{\alpha_1\over \alpha_2}x^2 
	 \over \ddots}
	 }}
	 }}
$$
where 
$ v_{2k+1}=\alpha_{k+1}/\alpha_k $,  
$ v_{2k+2}=\alpha_{k}/\alpha_{k+1} $
and $\alpha_n$ is defined by
$$
\sum_n \alpha_n x^n = {1+zx \over 1+yx + x^2}.
$$
Thus the Hankel determinants are
$$H(g)=(\alpha_0, \alpha_0, \alpha_1, \alpha_1, \alpha_2, \alpha_2, \alpha_3, \alpha_3, \ldots).$$

The proof of Theorem \ThmJFGene\ is based on the following generalization. 


\proclaim Theorem \ThmJFGeneLevel. 
Let $a_n, b_n, d_n, \alpha_n$ be numbers defined by
the following generating functions
{
\def\deno{(1-x)(1-(y^2-2)x+x^2)}
$$
\leqalignno{
	\sum_{n\geq 0} a_n x^n &= {(1-yz+z^2)(1-(y^2-2)x+x^2)    \over \deno},\cr
	\sum_{n\geq 0} b_n x^n &= {-z+y(yz-1)x-z(yz-1)x^2 \over \deno},\cr
	\sum_{n\geq 0} d_n x^n &= {-1-(1+z^2-2yz)x-z^2x^2  \over \deno},\cr
	\sum_n \alpha_n x^n  &= {1+zx \over 1+yx + x^2}\cr
}
$$
}
and $f_n(x)$ by
$$
(a_nx^4 -sb_nx^3+ b_n x^2) f_n(x)^2 + ((yd_n-2b_n) x^2 -sd_nx+ d_n) f_n(x) =d_n.
\leqno{(\EqJFQuadra)}
$$
Then
$$
f_n(x)=\displaystyle{1\over 1- sx -
     \displaystyle{ 
			 {\alpha_{n+1}\over \alpha_n} x^2\over 
	 1-  
	\displaystyle{
	{\alpha_n\over \alpha_{n+1}}x^2
		\over 1- sx -
     \displaystyle{ 
			 {\alpha_{n+2}\over \alpha_{n+1}} x^2\over 
			 1-{{\alpha_{n+1}\over \alpha_{n+2}}x^2 
	 \over \ddots}
	 }}
	 }}
$$
\medskip

When $n=0$, we have $a(0)=1-yz+z^2;\; b(0)=-z;\; d(0)=-1$.
Solving {(\EqJFQuadra)} yields Theorem \ThmJFGene.
For proving Theorem \ThmJFGeneLevel, we first convert it to 
Theorem \ThmJFGeneLevelK, in which the coefficients
are given by explicit formulas. This conversion is done by the
following change of variables:
$$
t+t^{-1}=y; \quad
t={y+\sqrt{y^2-4}\over 2}; \quad
K=(-t)^n.
\leqno{(\EqJFChangeVar)}
$$

\proclaim Theorem \ThmJFGeneLevelK.
Let $a(K), b(K), d(K), \alpha(K)$ be numbers defined by
$$
\leqalignno{
	a(K)&=
	{{{z}^{2}-(t+t^{-1})z+1}}
\cr
b(K)&=
-{\frac { \left( {t}^{2}+1 \right)  \gamma_2\gamma_3}{ \gamma_1^{2} 
}}-
	{\frac {{t} \gamma_2^{2}}{K^2 \gamma_1^{2} }}
	-{\frac {K^2 \gamma_3^{2}t}{ \gamma_1^{2} }}
\cr
d(K)&=
-\,{\frac {2t \gamma_2\gamma_3 }{ \gamma_1^{2} }}
	-{\frac { \gamma_2^2}{K^2 \gamma_1^{2} }}
	-{\frac {K^2{t}^{2} \gamma_3^{2}}{ \gamma_1^{2} }}\cr
\alpha(K)&=
-{\frac { \gamma_2 }{(1-{t}^{2})K}}  - {\frac {K \gamma_3 t}{1-{t}^{2}}}
\cr
}
$$
where $\gamma_1=(t-1)(t+1);\; \gamma_2=zt-1; \; \gamma_3=t-z$,
and $f(x; K)$ by
$$
\leqalignno{
&\qquad \qquad \bigl(a(K)x^4 -sb(K)x^3+ b(K) x^2\bigr) f(x;K)^2 
	  &(\EqJFChangeVarK)  \cr
&  + \bigl(((t+t^{-1})d(K)-2b(K)) x^2 -sd(K)x+ d(K)\bigr) f(x;K) =d(K). \cr
}
$$
Then
$$
f(x;K)=\displaystyle{1\over 1- sx-
     \displaystyle{ 
	 {\alpha(-tK)\over \alpha(K)} x^2\over 
	 1-  
	\displaystyle{
		{\alpha(K)\over \alpha(-tK)}x^2
		\over 1- sx-
     \displaystyle{ 
	 {\alpha(t^2K)\over \alpha(-tK)} x^2\over 
	 1-{{\alpha(-tK)\over \alpha(t^2K)}x^2 
	 \over \ddots}
	 }}
	 }}
$$


{\it Proof}.
Solving (\EqJFChangeVarK)  yields 
$$
f(x;K)=
{\frac { \left( (sx-1)t{ Q_1}- \left( t-1 \right)  \left( t+1 \right) { Q_2}\,{x}^{2}+{ Q_1}\,\sqrt {{ Q_0}}
	\right) { Q_1}}{2x^2 \left( t(1-sx){ Q_3}\,{ Q_1}+{K}^{2} \left( t-1 \right) ^{2} \left( t+1 \right) ^{2} \left( zt-1 \right) 
					  \left( t-z \right) {x}^{2} \right) }}
$$
where
$$
\leqalignno{
	Q_0&= \left( sxt-t-{x}^{2}-2\,{x}^{2}t-{x}^{2}{t}^{2} \right)  \left( sxt-t-{x}^{2}+2\,{x}^{2}t-{x}^{2}{t}^{2} \right); \cr
Q_1&= -1+zt+{K}^{2}{t}^{2}-{K}^{2}zt;\cr
Q_2&= {K}^{2}{t}^{2}-zt-{K}^{2}zt+1;\cr
Q_3&= -{K}^{2}z+z{t}^{2}-t+{K}^{2}t.\cr
}
$$
Then, we can verify 
$$
f(x;K)=\displaystyle{1\over 1- sx-
     \displaystyle{ 
	 {\alpha(-tK)\over \alpha(K)} x^2\over 
 1-{\alpha(-tK)\over \alpha(K)}x^2  f(x; -tK) }}.\qed
$$

\proclaim Proposition 3.6.
Let
$$
f=f(x)={1\over x^4}\sum_{k=1}^\infty {x^{2^{k+1}} \over 1-x^{2^k}}.
$$
Then, $H(f) \equiv (1,1,1,1,1,1,0,0)^* \pmod 2$.

\goodbreak
{\it Proof}. 
We successively have 
$$
\displaylines{
x^8f(x^2)=\sum_{k=2}^\infty {x^{2^{k+1}} \over 1-x^{2^k}}
=\sum_{k=1}^\infty {x^{2^{k+1}} \over 1-x^{2^k}} - {x^4\over 1-x^2}
=x^4f(x) - {x^4\over 1-x^2},\cr
x^4f(x^2)=f(x) - {1\over 1-x^2},\cr
x^4f(x)^2\equiv f(x) - {1\over 1-x^2} \pmod2,\cr
f(x)\equiv {1-\sqrt{1-{4x^4\over 1-x^2}}\over 2x^4} \pmod 2.\cr
}
$$
Let $g$ be the right-hand side of the above equation. By Proposition 3.1 we have
$$
H(g)=(1,1,1,1,-1-1,-2,-2,1,1,3,3,-1,-1,-4,-4, 1,1,5,5,\ldots).
$$
In other words,
$$
\leqalignno{
	H_{4k}(g)&= H_{4k+1}(g)=(-1)^k, \cr
	H_{4k+2}(g)&= H_{4k+3}(g)=(-1)^k(k+1), \cr
}
$$
so that
$$
H(f)\equiv H_n(g)\equiv (1,1,1,1,1, 1, 0, 0)^* \pmod 2\qed
$$

\proclaim Proposition 3.7.
Let
$$
f=f(x)={1\over x^4}\sum_{k=0}^\infty {x^{2^{k+2}} \over 1+x^{2^k}}.
$$
Then, $H(f) \equiv (1,1,0,0)^* \pmod 2$.

{\it Proof}.
Using the method described in the proof of Proposition 3.6, we derive 
$$
f(x)\equiv {1-\sqrt{1-{4x^4\over 1+x}}\over 2x^4} \pmod 2.
$$
Let $g$ be the right-hand side of the above equation. By Proposition 3.2 
$$
\displaylines{
H(g)=(1,1,0,0,-1,-1,0,0)^*.\cr
\noalign{\hbox{In other words,}}
H_{4k}(g)= H_{4k+1}(g)=(-1)^k, \cr
H_{4k+2}(g)= H_{4k+3}(g)=0. \cr
	\noalign{\hbox{Hence,}}
H(f)\equiv H_n(g)\equiv (1,1, 0, 0)^* \pmod 2.\qed\cr
}
$$

\proclaim Proposition 3.8. 
Let
$$
f=f(x)={1\over x^2}\sum_{k=0}^\infty {x^{2^{k+1}} \over 1+x^{2^{k+1}}}.
$$
Then $H(f) \equiv (1,1,0,0,1,1)^* \pmod 2$.

{\it Proof}. 
Using the method described in the proof of Proposition 3.6, we derive 
$$
f(x)\equiv {1-\sqrt{1-{4x^2\over 1+x^2}}\over 2x^2} \pmod 2.
$$
Let $g$ be the right-hand side of the above equality 
and let $\alpha_n$ be defined by
$$
\sum_n \alpha_n x^n = {1-x \over 1-x + x^2} = (1, 0, -1, -1, 0, 1)^*. 
$$
By Theorem 3.3 we have
$$H(g)=(1, 1, 0,0,  -1,-1, -1,  -1,0, 0,1, 1)^*. $$
Hence
$$H(f) \equiv H(g) \equiv (1,1,0,0,1,1)^* \pmod 2.\qed$$

\medskip


\section{4. One sequence, two modulos}   

Theorem 1.3 is proved in Section 2 by using the fractional congruence. In fact,
the fractional congruence can be avoided.
\smallskip
{\it Proof of Theorem 1.3}.
We successively have 
$$\displaylines{
	P_3(x)=(1-x)P_3(x^3)\equiv (1-x)P_3(x)^3 \pmod 3,\cr
P_3(x) ( 1- (1-x)P_3(x)^2) \equiv 0 \pmod 3,\cr
 1- (1-x)P_3(x)^2 \equiv 0 \pmod 3,\cr
P_3(x)^2 \equiv {1\over 1-x } \pmod 3,\cr
P_3(x)^2 \equiv {1\over (1-x)(1+3x) } \pmod 3,\cr
\Bigl(P_3(x) - \sqrt{1\over (1-x)(1+3x) }\Bigr) 
   \Bigl(P_3(x) + \sqrt{1\over (1-x)(1+3x) }\Bigr) \equiv 0 \pmod3.\cr 
}$$
We then have
$$P_3(x) \equiv \sqrt{1\over (1-x)(1+3x)} \pmod 3 \leqno{(4.1)} $$
by using the value of $P(0)$.
The right-hand side of the above equation has integral coefficients and its $J$-fraction is equal to:
$$
\sqrt{1\over (1-x)(1+3x)}=\JFrac{
2, (1)^*}{
	(-1)^*
}
,\leqno{(4.2)}
$$
so that $H_n(P_3)\equiv 2^{n-1} \equiv (-1)^{n-1} \pmod 3$ by (1.3). \qed

\medskip
Next we will prove the ``one sequence, two modulos" theorem 1.4.
We need the following lemma.

\proclaim Lemma 4.1.
We have
$$
\sqrt{1-4x} \equiv 1+2 \sum_{k=0}^\infty x^{2^k} \pmod 4.
$$

{\it Proof}.
The following expansion is well known (See [St99, WiCa] for example) 
$$
{ 1-\sqrt{1-4x} \over 2x} = \sum_{k=0}^\infty C_n x^n ,
$$
where $C_n={1\over n+1} {2n\choose n}$ is the Catalan number. 
It is easy to see that $C_n\equiv 1 \pmod 2$ if and only if $n=2^k-1$ for some integer $k$ [AK73], knowing, for instance, that $C_n$ is the number of binary trees with $n$ vertices. \qed

\medskip

{\it Proof of Theorem 1.4}.
We need to prove
$$
\leqalignno{
	\sqrt{1\over (1-x)(1+3x)} &\equiv \prod_{k=0}^\infty (1-x^{2^k})  \pmod 4,&(4.3)\cr
	\sqrt{1\over (1-x)(1+3x)} &\equiv \prod_{k=0}^\infty (1-x^{3^k})  \pmod 3.&(4.4)\cr
}
$$
The second equality is just relation (4.1).  
For proving the first equality
let $f(x)$ be the left-hand side of (4.3). By Lemma 4.1, we get
$$
\displaylines{
(1-x) f(x) = \sqrt{1-x\over 1+3x} = \sqrt{1-{4x\over 1+3x}}
\equiv 1+2 \sum_{k=0}^\infty \bigl({x\over 1+3x}\bigr)^{2^k} \pmod 4,\cr
(1-x) f(x) \equiv 1+2 \sum_{k=0}^\infty \bigl({x\over 1+x}\bigr)^{2^k} \pmod 4,\cr
}$$
and
$$
\leqalignno{
(1-x^2) f(x^2) 
&\equiv 
1+2 \sum_{k=0}^\infty \bigl({x^2\over 1+x^2}\bigr)^{2^k}  &(4.5)\cr
&\equiv 
1+2 \sum_{k=0}^\infty \bigl({x\over 1+x}\bigr)^{2^{k+1}} \qquad \hbox{[By Lemma 2.1]}  \cr
&\equiv
(1-x)f(x)-{2x\over 1+x} \pmod 4. \cr
}
$$
Let $P_2(x)$ be the right-hand side  of (4.3). Then, 
$$
\displaylines{
(1-x)P_2(x^2)=P_2(x),\cr
\noalign{\hbox{which implies on one hand $(1-x)(P_2(x))^2\equiv P_2(x)\pmod 2$, hence}}
P_2(x)\equiv {1\over 1+x} \pmod 2,\cr
\noalign{\hbox{and on theother hand}}
(1-x^2)P_2(x^2) = (1+x) (1-x) P_2(x^2) = (1+x)P_2(x).\cr
}
$$
Hence,
$$
(1-x^2)P_2(x^2) 
- (1-x)P_2(x) 
\equiv  {2x\over 1+x}
\pmod 4.
\leqno{(4.6)}
$$
Taking the difference of (4.6) and (4.5) yields
$$
(1-x^2)f(x^2)-(1-x^2)P_2(x^2) \equiv (1-x)(f(x)-P_2(x)) \pmod 4,
$$
and
$$
f(x)-P_2(x) \equiv (1+x)(f(x^2)-P_2(x^2))  \pmod 4.\leqno{(4.7)}
$$
By applying (4.7) recursively we get $f(x)-P_2(x) \equiv 0 \pmod 4$, 
since $f(0)=P_2(0)=1$.  \qed


\goodbreak
\bigskip
{\bf Acknowledgements}. The author should like to thank Zhi-Ying Wen who 
suggested that I study
the Hankel determinants of the Thue-Morse sequence back to 1991,
and who invited me to Tsinghua University where the paper was finalized.
The author also thanks 
Yann Bugeaud for reanimating the topic and valuable discussion.  

\bigskip

\bigskip
\centerline{\bf References}
\bigskip

{
\eightpoint

\article AK73|Alter, R., Kubota, K|Prime and prime power divisibility of Catalan numbers|J. Combin. Theory Ser. A|15|1973|243--256|

\article APWW|Allouche, J.-P.; Peyri\`ere, J.; Wen, Z.-X.; Wen, Z.-Y|
Hankel determinants of the Thue-Morse sequence|Ann. Inst. Fourier
$($Grenoble$)$|48|1998|1--27|

\divers BH13|Bugeaud, Yann; Han, Guo-Niu|A combinatorial proof of the non-vanishing of Hankel determinants of the Thue--Morse sequence. {\it Preprint}, 16 pages, {\oldstyle 2013}, {\tt  http://www-irma.u-strasbg.fr/\char126guoniu/papers/}|

\article Co13|Coons, Michael|On the rational approximation of the sum of the reciprocals of the Fermat numbers|Ramanujan J.|30|2013|39--65|

\article Fl80|Flajolet, Philippe|Combinatorial aspects of continued fractions|Discrete Math.|32|1980|125--161|

\divers Kr98|Krattenthaler, Christian|Advanced determinant calculus,
{\sl S\'em. Lothar. Combin.},
{\bf B42q} ({\oldstyle1998}), 67pp|

\article Kr05|Krattenthaler, Christian|Advanced determinant calculus: A complement|Linear Algebra and Appl.|411|2005|68--166|

\livre Mu23|Muir, T|The theory of
determinants in the historical order of development, {\rm
4~vols}|Macmillan, London, {\oldstyle  1906}--{\oldstyle 1923}|

\divers Ru06|Rubey, Martin|Extended Rate, more GFUN, {\sl Fourth Colloquium on Mathematics and Computer Science}, DMTCS Proc. AG. {\oldstyle 2006}, pp. 431--434|

\livre St99|Stanley, Richard P.|Enumerative Combinatorics, vol. 2|
Cambridge University Press, {\oldstyle 1999}|

\divers Vi83|Viennot, X|Une th\'eorie combinatoire des polyn\^omes 
orthogonaux g\'en\'eraux, {\sl UQAM, Montreal, Quebec}, {\oldstyle 1983}| 

\livre Wa48|Wall, H. S|Analytic theory of continued fractions|Chelsea Publishing Company, Bronx, N.Y., {\oldstyle 1948}|

\divers WiCa|Wikipedia|Catalan number, {\it revision September 2, {\oldstyle 2013}}|

\bigskip

}
\bigskip\bigskip
\hbox{
\vtop{\halign{#\hfil\cr
IRMA, UMR 7501\cr
Universit\'e de Strasbourg et CNRS\cr
7 rue Ren\'e Descartes\cr
67084 Strasbourg, France\cr
\noalign{\smallskip}
{\tt guoniu.han@unistra.fr}\cr}}}


\end